\documentclass[12pt]{article}

\newcommand{\be}{\begin{equation}}
\newcommand{\ee}{\end{equation}}
\newcommand{\proof}{\noindent {\bf Proof.\ \ }}

\newtheorem{theorem}{Theorem}

\newtheorem{corollary}{Corollary}

\mathchardef\inn="3232
\renewcommand{\in}{\mbox{$\,\inn\,$}}

\begin{document}


\title{\bf Upper Bounds on the Automorphism Group of a Graph \\
Discrete Mathematics 256 (2002) 489-493.
}
\author{   {\bf Ilia Krasikov}\\
	    Brunel University\\
    Department of Mathematical Sciences\\
Uxbridge UB8 3PH United Kingdom\\
   {\bf Arie Lev }\\
  The Academic College of Tel-Aviv-Yaffo,\\ Tel-Aviv Israel 64044\\
   {\bf Bhalchandra D. Thatte}\\
   450 N Mathilda Av., M107, Sunnyvale CA 94085 USA }

\date{}


\maketitle
\vspace*{4ex}

\begin{center}
{\bf Abstract}

We give upper bounds on the order of the automorphism group of a
simple graph

\end{center}

\baselineskip3ex

\vspace{16ex}

\thispagestyle{empty}
\addtocounter{page}{-1}

\baselineskip3.3ex

%
%
%
%
%
%
%

In this note we present some upper bounds on the order of the
automorphism group of a graph, which is assumed to be simple,
having no loops or multiple edges. Somewhat surprisingly, we did
not find such bounds in the literature and the goal of this paper
is to fill this gap. As a matter of fact, implicitly such bounds
were contained in works dealing with the edge reconstruction
conjecture and are the corollaries of a simple theorem which is
presented below (Theorem 1). Therefore we bring together a few
results spread in different, sometimes in difficult to reach,
sources (see Theorem 2 below). In Theorem 3 we derive a new bound,
based on the notion of a \it greedy spanning tree \rm . This new
bound improves, in many cases, the bounds (1) and (2) of Theorem
2.
\par
We will use the following notation. Let $F$ be a spanning subgraph
of a fixed copy of a graph $G$. The number of embeddings of $F$ in
$G$, that is the number of labeled copies of $F$ in $G$, is
denoted by $|F \rightarrow G|$. Clearly $|F \rightarrow G|=s(F
\rightarrow G)aut( F),$ where $s(F \rightarrow G)$ is the number
of subgraphs of $G$ isomorphic to $F$ and $aut(F)$ is the order of
the automorphism group of $F$. We also use $n=n(G)$ for the number
of vertices and $e=e(G)$ for the number of edges of $G$. As usual,
$\Delta_G$, $\delta_G$ and $d_G$ stand for the maximum, the
minimum and the average degree of $G$ respectively. The degree of
a vertex $v \in G$ is denoted by $d_G (v)$.
\begin{theorem}
Let $F$ be a spanning subgraph of a graph $G$, Then $$aut(G) \le
|F \rightarrow G|=s(F \rightarrow G)aut(F).$$
\end{theorem}
\proof Let $\phi : G \rightarrow G$ be an automorphism of $G$ and
let $F_1$ be a fixed copy of $F$ in $G$. Then, as $F$ is a
spanning subgraph of $G$, $\phi$ is completely determined by the
knowledge of $\phi(F_1)$. Since the number of different images
$\phi(F_1)$ does not exceed $|F \rightarrow G|$, the result
follows.

Some relevant estimates of $|F \rightarrow G|$, $ s(F \rightarrow G)$ 
and $aut(F)$
for graphs in general and for special families of graphs
are known and have been obtained mainly in connection with
the edge reconstruction conjecture. We try to collect them in the 
following
\begin{theorem}
Let $G$ be a connected graph, then
\be
\label{u1}
aut(G) \le n (\Delta_G)! \, ( \Delta_G -1)^{n- \Delta_G-1}
\ee

Let $T$ be a spanning tree in $G$, then
\be
\label{u2} aut(G) \le \frac{\Delta_T}{\Delta_G} \, (d_G)^n
\prod_{v \in V(G)} (d_T(v)-1)! \ee

Let $p=p(G)$ be the path covering number of a graph, i.e. the minimum 
number
of vertex-disjoint paths containing all vertices of $G$. Then
\be
\label{u7}
aut(G) \le 2p \, n^{2p} (2^{7/8} 6^{1/24})^{e-n}
\ee
\be
\label{u8}
aut(G) \le (d_G)^n 
((\Delta_G-1)!)^{\frac{e-n+3-2\delta_G}{(\delta_G-1)(\Delta_G-2)} } ,
\ee
provided $\delta_G \ge 2$, $\Delta_G \ge 3$.

Let $G$ be either a square of a graph or a three-connected planar 
graph, then
\be
\label{u3}
aut(G) \le 3 \, \frac{  2^{\frac{n-2}{2}} (d_G)^n}{\Delta_G }
\ee

Let $G$ be a $K_{1,m}$-free graph, then
\be
\label{u4}
aut(G) \le \frac{ (m-1)! ((m-2)!)^{\frac{n}{m-2}} \, (d_G)^n}{\Delta_G}
\ee

If $G$ has a hamiltonian path then
\be
\label{u5}
aut(G) \le n (\frac{e}{n-1})^{n-1}
\ee

\be
\label{u6}
aut(G) \le 2 n^2 (2^{7/8} 6^{1/24})^{e-n}
\ee
\end{theorem}
\proof
Everywhere in the sequel $T$ is a spanning tree in $G$.
The bound (\ref{u1}) is just Caunter and Nash-Williams' estimate for 
$|T \rightarrow G|$,
see \cite{caun} and \cite{bondyp, fio}.\\
It has been shown in \cite{ckr1} that
\be
\label{fs}
s(T \rightarrow G) \le \frac{\prod_{v \in G} d_G (v)}{\Delta_G} \le 
\frac{d_G^n}{\Delta_G},
\ee
and
\be
\label{fa} aut(T)\leq \Delta_T \, \prod_{v \in T} (d_T(v)-1)!, \ee
giving (\ref{u2}), see also \cite{bondy}.\\ If $G$ satisfies
$\delta_G \ge 2$, $\Delta_G \ge 3$, than there is a spanning tree
$T$ in $G$ such that \cite{ckr1} $$aut(T) \le \Delta _G \,
((\Delta_G-1)!)^{\frac{e-n+3-2\delta_G}{(\delta_G-1)(\Delta_G-2)}
},$$ This gives (\ref{u8}) by (\ref{fs}).\\ Concerning (\ref{u3})
notice that in both cases the corresponding graphs have a spanning
tree of maximum degree at most 3. For the square of a graph this
has been proved in \cite{ckr2} and for three-connected planar
graphs this is a classical result of Barnette \cite{bar}. These
yield (\ref{u3}) by (\ref{fs}) and (\ref{fa}) since the maximum of
the product in (\ref{fa}) is attained then the tree has the
maximal possible number $\frac{n-2}{2}$ of vertices of degree
$3$.\\ The required estimates for $K_{1,m}$-free graphs giving
(\ref{u4}) has been established in \cite{ckr1}. Namely, a
$K_{1,m}$-free graph has a spanning tree of maximum degree at most
$m$. Moreover, such a tree can be modified to have $aut(T) \le
(m-1)! ((m-2)!)^{\frac{n}{m-2}}.$\\ The inequality (\ref{u7}) due
to Pyber \cite{pyb}. If $P$ is a hamiltonian path Lov\'{a}sz
proved \cite{lov} (see also \cite{bondyp}) $s(P \rightarrow G) \le
\frac{n}{2} (\frac{e}{n-1})^{n-1}.$ Since $aut(P)=2$ this yields
(\ref{u5}). Finally (\ref{u6}) follows from (\ref{u7}) with $p=1$.

We derive now another bound on $aut(G)$ (see Theorem 3 below).
First, we shall define the notion of a \it greedy spanning tree,
\rm $T(v_0,v_1,\dots v_s)$, of a connected graph $G$ by the
following construction:
\par
We shall define the sequence of vertices $v_0,v_1,v_2,\dots ,v_s$
of $G$ and the corresponding sequence $T_0,T_1,\dots T_s$ of trees
as follows: Let $v_0$ be any vertex of $G$ and let $T_0$ be the
tree containing $v_0$ and all the edges of $G$ which are adjacent
to $v_0$ (we mean that if a subgraph contains an edge, then it
contains also its end vertices). Note that $T_0$ is actually a
star with central vertex $v_0$. In order to construct $T_1$ choose
any leaf $v_1$ of $T_0$ having at least one adjacent edge which is
not adjacent to any vertex of $T_0=\{ v_1\} $, and add to $T_0$
all the edges adjacent to $v_1$ which are not adjacent to any
vertex in $V(T_0)-\{ v_1\} $. Denote the resulting tree by $T_1$.
Continue this construction inductively: given $T_{i-1}$, let $v_i$
be a leaf of $T_{i-1}$ having an adjacent edge which is not
adjacent to any vertex of $T_{i-1}-\{ v_i\} $, and add to
$T_{i-1}$ all the edges which are adjacent to $v_i$ and which are
not adjacent to any vertex of $V(T_{i-1})-\{ v_i\} $. Denote the
resulting graph by $T_i$. This construction is completed at step
$s$, when for every leaf $v$ of $T_s$, each edge of $G$ which is
adjacent to $v$, is also adjacent to a vertex in $V(T_s)-\{ v\} $.
\par
It is easy to see that for a connected graph $G$, the above
(greedy) construction results in a spanning tree $T_s$ of $G$.
This spanning tree will be called a \it greedy spanning tree \rm
of $G$ and denoted by $T=T(v_0,v_1,\dots ,v_s)$, where
$v_0,v_1,\dots ,v_s$ is the sequence of vertices used in the above
construction of $T$.
\par
Using the above notation, we derive the following bound for
$autG$.

\

\begin{theorem}
Let $G$ be a connected simple graph with $n$ vertices and let
$T=T(v_0,v_1,\dots ,v_s)$ be a greedy spanning tree of $G$. Denote
by $n_1$ the length of the orbit of $v_0$ under the action of the
automorphism group of $G$. Then $$autG\leq
n_1(d(v_0))!\prod_{i=1}^s(d_T(v_i)-1)!$$ In particular, for any
greedy spanning tree $T$ of $G$ we have: $$autG\leq
n(d(v_0))!\prod_{v\in V(G)}(d_T(v)-1)!$$
\end{theorem}

\

\proof Let $\Gamma $ be the automorphism group of $G$. Given
vertices $u_1,u_2,\dots ,u_r$ of $G$, denote by $C_{\Gamma
}(u_1,u_2,\dots ,u_r)$ the subgroup of $\Gamma $ which fixes
$u_1,u_2,\dots ,u_r$. Then we have $aut(G)=n_1|C_{\Gamma }(v_0)|$.
Since $C_{\Gamma }(v_0)$ acts on the set $N(v_0)$ (the set of all
neighbors of $v_0$ in $G$), and since $v_1$ is a neighbor of
$v_0$, we have $|C_{\Gamma }(v_0)|\leq d(v_0)|C_{\Gamma
}(v_0,v_1)|$ (equality holds if and only if $C_{\Gamma }(v_0)$ is
transitive on $N(v_0)$). Denote $N(v_0)=\{ v_1,u_2,\dots
u_{d(v_0)}\} $. Then, we have: $|C_{\Gamma }(v_0,v_1)|\leq
(d(v_0)-1)|C_{\Gamma }(v_0,v_1,u_2)|\leq
(d(v_0-1))(d(v_0-2))|C_{\Gamma }(v_0,v_1,u_2.u_3)|\leq \dots \leq
(d(v_0)-1)!|C_{\Gamma }(v_0,v_1,u_2,\dots ,u_{d(v_0)})|$. Whence
$aut(G)\leq n_1(d(v_0))!|C_{\Gamma }(\{ v_0\} \cup N(v_0))|$.
\par
Since $v_1$ is adjacent to $v_0$ in $T$, we have that $C_{\Gamma
}(\{ v_0\} \cup N(v_0))$ acts on $N_T(v_1)-\{ v_0\} $. It follows
by the arguments used in the preceding paragraph that $|C_{\Gamma
}(\{ v_0\} \cup N(v_0))|\leq (d_T(v_1)-1)!|C_{\Gamma }(
N(v_0)\cup N_T(v_1))|$, and consequently, $aut(G)\leq
n_1(d(v_0))!(d_T(v_1)-1)!|C_{\Gamma }(N(v_0)\cup N_T(v_1))|$. The
theorem now follows by repeating the above arguments for the
vertices $v_2,v_3,\dots ,v_s$.

The following corollary is a straightforward result of Theorem 3.

\begin{corollary}
Denote $r=\lfloor \frac{n-\Delta_G-1}{\Delta_G-1}\rfloor $,
$\alpha =n-r(\Delta_G-1)$ (clearly $\alpha < \Delta_G -1$). Then
$$autG\leq n\alpha !\Delta_G![(\Delta_G-1)!]^r.$$
\end{corollary}

\par
\it Remark \rm It is easily verified that $aut(K_n)=n!$,
$aut(K_{m,m})=2(m!)^2$ and $aut(K_{p,q})=p!q!$ for $p\neq q$.
Applying Theorem 3 for these graphs, we have that the bound of
Theorem 3 is exact (i.e., the corresponding inequality is actually
an equality). On the other hand, except for the case of formula
(1) applied for $aut(K_n)$, the inequalities of Theorem 2 are not
exact in the above cases.


\newpage

\begin{thebibliography}{8}

\bibitem{bar}
D.\,Barnette, {\em Trees in polytopal graphs},
Canad. J. Math. 18, 1966, pp. 731-736.

\bibitem{bondy}
J.\,A.\,Bondy, {\em A graph reconstructor's manual},
in Surveys in Combinatorics, 1991, London Math. Soc. Lecture Note Ser. 
166,
Cambridge Univ. Press. Cambridge 1991, pp. 221-252.

\bibitem{bondyp}
J.\,A.\,Bondy, {\em The reconstruction of graphs}, preprint,
Dept. of Combinatorics and Optimization, University of Waterloo, 1983.

\bibitem{caun} J.\, Caunter and C.\,St.\,J.\,A.\,Nash-Williams, {\em
Degree conditions for edge reconstruction}, preprint, 1982.

\bibitem{ckr1}
Y.\,Caro, I.\,Krasikov, Y.\,Roditty, {\em  Spanning trees and some 
edge-reconstructible graphs},
Ars Combinatoria, 20-A, 1985, pp. 109-118.

\bibitem{ckr2}
Y.\,Caro, I.\,Krasikov, Y.\,Roditty, {\em On the largest subtree
of a given maximum degree in connected graphs},  J. Graph Theory
15, 1991, pp. 7-13.

\bibitem{el}
M.\,N.\,Ellingham, {\em Recent progress in edge reconstruction},
Congressus Numerantium 62, 1988, pp.3-20.

\bibitem{fio}
S.\,Fiorini and J.\,Lauri, {\em Edge reconstruction of graphs with 
topological properties},
in Combinatorial Mathematics (Marseille-Luminy, 1981),
North-Holland Math. Stud. 75, North-Holland, Amsterdam-New York, 1983, 
pp. 285-288.

\bibitem{k1}
I.\,Krasikov, {\em A note on the edge-reconstruction of $K_{1,m}$-free 
graphs},
J. Comb. Theory (B), Vol.49, No 2, 1990, pp. 295-298.

\bibitem{kr1}
I.\,Krasikov, Y.\,Roditty, {\em  Recent applications of Nash-Williams
lemma to the edge-reconstruction conjecture},
Ars Combinatoria, 29-A, 1990, pp. 215-224.

\bibitem{lov}
L.\,Lov\'{a}sz, Some problems of graph theory, Matematikus Kurir, 1983.

\bibitem{pyb}
L.\,Pyber, {\em The edge reconstruction of hamiltonian graphs},
J. Graph Theory 14, 1990, pp.173-179.
\end{thebibliography}
\end{document}